\documentclass[reqno,onecolumn,12pt]{amsart}
\usepackage{amsfonts}
\usepackage{amssymb}
\usepackage{amsmath}
\usepackage{graphicx}
\usepackage{amscd}

\newtheorem{theorem}{Theorem}
\theoremstyle{plain}

\newtheorem{lemma}{Lemma}

\newtheorem{proposition}{Proposition}
\newtheorem{remark}{Remark}

\numberwithin{equation}{section}

\input{tcilatex}

\begin{document}
\author{}
\title{}
\maketitle

\begin{center}
\bigskip \thispagestyle{empty}\textbf{EXTENDED FERMIONIC }$p$\textbf{-ADIC }$%
q$-\textbf{INTEGRALS ON }$%
\mathbb{Z}
_{p}$\textbf{\ IN CONNECTION WITH APPLICATIONS OF UMBRAL CALCULUS}

\bigskip \textbf{By}

\hspace{0.05cm}

$^{\dag }$\textbf{Serkan ARACI, }$^{\dag }$\textbf{Mehmet ACIKGOZ and }$%
^{\ddag }$\textbf{Erdo\u{g}an \c{S}EN}

\hspace{-0.5cm}

$^{\dag }$University of Gaziantep, Faculty of Arts and Science, Department
of Mathematics, 27310 Gaziantep, Turkey

\hspace{0.05cm}

$^{\ddag }$Department of Mathematics, Faculty of Science and Letters, Namik
Kemal University, 59030 Tekirda\u{g}, Turkey

\hspace{0.05cm}

mtsrkn@hotmail.com; acikgoz@gantep.edu.tr; erdogan.math@gmail.com

\hspace{-0.5cm}

\textbf{{\Large {Abstract}}}

\hspace{-0.5cm}
\end{center}

The purpose of this paper is to derive some applications of umbral calculus
by using extended fermionic $p$-adic $q$-integral on $%
\mathbb{Z}
_{p}$. From those applications, we derive some new interesting properties on
the new family of Euler numbers and polynomials. That is, a systemic study
of the class of Sheffer sequences in connection with generating function of
the weighted $q$-Euler polynomials are given in the present paper.

\bigskip


\hspace{-0.5cm}


\section{\textbf{Preliminaries}}

Suppose that $p$ be a fixed odd prime number. Throughout this work we use
the following notations, where $%
\mathbb{Z}
_{p}$ we denote the ring of $p$-adic rational integers, $%
\mathbb{Q}
$ denotes the field of rational numbers, $%
\mathbb{Q}
_{p}$ denotes the field of $p$-adic rational numbers, and $%
\mathbb{C}
_{p}$ denotes the completion of algebraic closure of $%
\mathbb{Q}
_{p}$. Let $%
\mathbb{N}
$ be the set of natural numbers and $%
\mathbb{N}
^{\ast }=%
\mathbb{N}
\cup \left\{ 0\right\} $. The $p$-adic absolute value is defined by $%
\left\vert p\right\vert _{p}=p^{-1}$. Also, we assume that $\left\vert
q-1\right\vert _{p}<1$ is an indeterminate. Let $UD\left( 
\mathbb{Z}
_{p}\right) $ be the space of uniformly differentiable functions on $%
\mathbb{Z}
_{p}$. For $f\in UD\left( 
\mathbb{Z}
_{p}\right) $, the fermionic $p$-adic $q$-integral on $%
\mathbb{Z}
_{p}$ is defined by T. Kim, as follows:%
\begin{equation}
I_{-q}\left( f\right) =\int_{%
\mathbb{Z}
_{p}}f\left( \xi \right) d\mu _{-q}\left( \xi \right) =\lim_{n\rightarrow
\infty }\frac{1}{\left[ p^{n}\right] _{-q}}\sum_{\xi =0}^{p^{n}-1}\left(
-1\right) ^{\xi }f\left( \xi \right) q^{\xi }\text{.}  \label{euation 34}
\end{equation}%
where $\left[ x\right] _{q}$ is $q$-analogue of $x$ defined by%
\begin{equation*}
\left[ x\right] _{q}=\frac{q^{x}-1}{q-1}\text{ \textit{and} }\left[ x\right]
_{-q}=\frac{1-\left( -q\right) ^{x}}{1+q}\text{.}
\end{equation*}

We want to note that $\lim_{q\rightarrow 1}\left[ x\right] _{q}=x$ (for
details, see [1-31]).

By (\ref{euation 34}), we have%
\begin{equation}
qI_{-q}\left( f_{1}\right) +I_{-q}\left( f\right) =\left[ 2\right]
_{q}f\left( 0\right)  \label{equation 4}
\end{equation}%
where $f_{1}\left( \xi \right) :=f\left( \xi +1\right) $ (for details, see 
\cite{Kim 3}, \cite{Kim 8}).

Let us consider Kim's fermionic $p$-adic $q$-integral on $%
\mathbb{Z}
_{p}$ in the following form: for $\left\vert 1-\zeta \right\vert _{p}<1$%
\begin{equation}
I_{-q}^{\zeta }\left( f\right) =\int_{%
\mathbb{Z}
_{p}}\zeta ^{\xi }f\left( \xi \right) d\mu _{-q}\left( \xi \right)
=\lim_{n\rightarrow \infty }\frac{1}{\left[ p^{n}\right] _{-q}}\sum_{\xi
=0}^{p^{n}-1}\zeta ^{\xi }f\left( \xi \right) \left( -1\right) ^{\xi }q^{\xi
}\text{,}  \label{equation 47}
\end{equation}%
where $I_{-1}^{\zeta }\left( f\right) $ are called extended fermionic $p$%
-adic $q$-integral on $%
\mathbb{Z}
_{p}$.

Let us now consider $f_{1}\left( \xi \right) :=f\left( \xi +1\right) $, then
we develop as follows:%
\begin{eqnarray*}
-q\zeta I_{-q}^{\zeta }\left( f_{1}\right) &=&\lim_{n\rightarrow \infty }%
\frac{1}{\left[ p^{n}\right] _{-q}}\sum_{\xi =0}^{p^{n}-1}\zeta ^{\xi
+1}\left( -1\right) ^{\xi +1}f\left( \xi +1\right) q^{\xi +1} \\
&=&I_{-q}^{\zeta }\left( f\right) +\frac{\left[ 2\right] _{q}}{2}%
\lim_{n\rightarrow \infty }\left( -f\left( 0\right) -\xi ^{p^{n}}f\left(
p^{n}\right) q^{p^{n}}\right) \\
&=&I_{-q}^{\zeta }\left( f\right) -\left[ 2\right] _{q}f\left( 0\right) 
\text{.}
\end{eqnarray*}

Therefore, we have the following lemma.

\begin{lemma}
\label{Lemma 1}For $f\in UD\left( 
\mathbb{Z}
_{p}\right) $, we get%
\begin{equation*}
q\zeta I_{-q}^{\zeta }\left( f_{1}\right) +I_{-q}^{\zeta }\left( f\right) = 
\left[ 2\right] _{q}f\left( 0\right) \text{.}
\end{equation*}
\end{lemma}

Taking $f\left( \xi \right) =e^{t\left( x+\xi \right) }\in UD(%
\mathbb{Z}
_{p})$ in Lemma \ref{Lemma 1}, then we introduce the following expression:%
\begin{equation}
\int_{%
\mathbb{Z}
_{p}}\zeta ^{\xi }e^{t\left( x+\xi \right) }d\mu _{-q}\left( \xi \right) =%
\frac{\left[ 2\right] _{q}}{q\zeta e^{t}+1}e^{tx}=\sum_{n=0}^{\infty
}E_{n,\zeta }^{q}\left( x\right) \frac{t^{n}}{n!}\text{,}
\label{equation 35}
\end{equation}%
where we call $E_{n,\zeta }^{q}\left( x\right) $ as weighted $q$-Euler
polynomials. In the special case, $x=0$, $E_{n,\zeta }^{q}\left( 0\right)
:=E_{n,\zeta }^{q}$ are called weighted $q$-Euler numbers and the relation
between weighted $q$-Euler numbers and weighted $q$-Euler polynomials are
given by%
\begin{equation}
E_{n,\zeta }^{q}\left( x\right) =\sum_{l=0}^{n}\binom{n}{l}x^{l}E_{n-l,\zeta
}^{q}=\left( x+E_{\zeta }^{q}\right) ^{n}\text{,}  \label{equation 36}
\end{equation}%
with the usual of replacing $\left( E_{\zeta }^{q}\right) ^{n}$ by $%
E_{n,\zeta }^{q}$ is used. By (\ref{equation 35}), we note that 
\begin{equation}
E_{n,\zeta }^{q}=\int_{%
\mathbb{Z}
_{p}}\zeta ^{\xi }\xi ^{n}d\mu _{-q}\left( \xi \right) \text{ \textit{and} }%
E_{n,\zeta }^{q}\left( x\right) =\int_{%
\mathbb{Z}
_{p}}\zeta ^{\xi }\left( x+\xi \right) ^{n}d\mu _{-q}\left( \xi \right) 
\text{.}  \label{euation 11}
\end{equation}

By (\ref{equation 35}), we have%
\begin{equation}
E_{n,\zeta }^{q}\left( x\right) =\left[ 2\right] _{q}\sum_{m=0}^{\infty
}\left( -1\right) ^{m}q^{m}\zeta ^{m}\left( m+x\right) ^{n}\text{, \textit{%
for} }n\in 
\mathbb{N}
^{\ast }\text{.}  \label{euation 50}
\end{equation}

From this, we can define weighted $q$-Zeta function as follows:%
\begin{equation}
\lambda \left( s,x:q:\zeta \right) =\left[ 2\right] _{q}\sum_{m=0}^{\infty }%
\frac{\left( -1\right) ^{m}q^{m}\zeta ^{m}}{\left( m+x\right) ^{s}}
\label{euation 51}
\end{equation}

By (\ref{euation 50}) and (\ref{euation 51}), we derive the following
equation (\ref{euation 52}):%
\begin{equation}
\lambda \left( -n,x:q:\zeta \right) =E_{n,\zeta }^{q}\left( x\right) ,\text{ 
\textit{for any} }n\in 
\mathbb{N}
^{\ast }\text{.}  \label{euation 52}
\end{equation}

When we set as $q=\zeta =1$ in (\ref{euation 52}) which reduces to%
\begin{equation*}
\zeta _{E}\left( -n,x\right) =E_{n}\left( x\right)
\end{equation*}%
which is well known in \cite{Kim 14}.

By (\ref{equation 47}) and (\ref{equation 35}), we compute%
\begin{eqnarray*}
\int_{%
\mathbb{Z}
_{p}}\zeta ^{\xi }\left( x+\xi \right) ^{n}d\mu _{-q}\left( \xi \right)
&=&\lim_{m\rightarrow \infty }\frac{1}{\left[ dp^{m}\right] _{-q}}\sum_{\xi
=0}^{dp^{m}-1}\left( -1\right) ^{\xi }\zeta ^{\xi }\left( x+\xi \right)
^{n}q^{\xi } \\
&=&\frac{d^{n}}{\left[ d\right] _{-q}}\sum_{j=0}^{d-1}\left( -1\right)
^{j}\zeta ^{j}q^{j}\left( \lim_{m\rightarrow \infty }\frac{1}{\left[ p^{m}%
\right] _{\left( -q\right) ^{d}}}\sum_{\xi =0}^{p^{m}-1}\left( -1\right)
^{\xi }\left( \zeta ^{d}\right) ^{\xi }\left( q^{d}\right) ^{\xi }\left( 
\frac{x+j}{d}+\xi \right) ^{n}\right) \\
&=&\frac{d^{n}}{\left[ d\right] _{-q}}\sum_{j=0}^{d-1}\left( -1\right)
^{j}\zeta ^{j}q^{j}\int_{%
\mathbb{Z}
_{p}}\zeta ^{d\xi }\left( \frac{x+j}{d}+\xi \right) ^{n}d\mu _{-q^{d}}\left(
\xi \right) \text{,}
\end{eqnarray*}%
where $d$ is an odd natural number. So from the above%
\begin{equation}
\int_{%
\mathbb{Z}
_{p}}\zeta ^{\xi }\left( x+\xi \right) ^{n}d\mu _{-q}\left( \xi \right) =%
\frac{d^{n}}{\left[ d\right] _{-q}}\sum_{j=0}^{d-1}\left( -1\right)
^{j}\zeta ^{j}q^{j}\int_{%
\mathbb{Z}
_{p}}\zeta ^{d\xi }\left( \frac{x+j}{d}+\xi \right) ^{n}d\mu _{-q^{d}}\left(
\xi \right) \text{.}  \label{euation 53}
\end{equation}

By (\ref{euation 11}) and (\ref{euation 53}), we get%
\begin{equation}
E_{n,\zeta }^{q}\left( dx\right) =\frac{d^{n}}{\left[ d\right] _{-q}}%
\sum_{j=0}^{d-1}\left( -1\right) ^{j}\zeta ^{j}q^{j}E_{n,\zeta
^{d}}^{q^{d}}\left( x+\frac{j}{d}\right) \text{,}  \label{euation 54}
\end{equation}%
which plays an important role for studying regarding Measure theory on $p$%
-adic analysis.

Let we use the following notations, where $%
\mathbb{C}
$ denotes the set of complex numbers, $\mathcal{F}$ denotes the set of all
formal power series in the variable $t$ over $%
\mathbb{C}
$ with 
\begin{equation*}
\mathcal{F}=\left\{ f\left( t\right) =\sum_{k=0}^{\infty }a_{k}\frac{t^{k}}{%
k!}\mid a_{k}\in 
\mathbb{C}
\right\} ,
\end{equation*}%
$\mathcal{P}=%
\mathbb{C}
\left[ x\right] $ and $\mathcal{P}^{\ast }$ denotes the vector space of all
linear functional on $\mathcal{P}$, $\left\langle L\mid p\left( x\right)
\right\rangle $ denotes the action of the linear functional $L$ on the
polynomial $p\left( x\right) $, and it is well-known that the vector space
operation on $\mathcal{P}^{\ast }$ is defined by 
\begin{eqnarray*}
\left\langle L+M\mid p\left( x\right) \right\rangle &=&\left\langle L\mid
p\left( x\right) \right\rangle +\left\langle M\mid p\left( x\right)
\right\rangle , \\
\left\langle cL\mid p\left( x\right) \right\rangle &=&c\left\langle L\mid
p\left( x\right) \right\rangle ,
\end{eqnarray*}%
where $c$ is any constant in $%
\mathbb{C}
$ (for details, see \cite{Kim 2}, \cite{Dere 1}, \cite{Kim 11}, \cite{Kim 12}%
, \cite{Roman}).

The formal power series are known by%
\begin{equation*}
f\left( t\right) =\sum_{k=0}^{\infty }a_{k}\frac{t^{k}}{k!}\in \mathcal{F}
\end{equation*}%
which describes a linear functional on $\mathcal{P}$ as $\left\langle
f\left( t\right) \mid x^{n}\right\rangle =a_{n}$ for all $n\geq 0$ (for
details, see \cite{Kim 2}, \cite{Dere 1}, \cite{Kim 11}, \cite{Kim 12}, \cite%
{Roman}). In addition to%
\begin{equation}
\left\langle t^{k}\mid x^{n}\right\rangle =n!\delta _{n,k},
\label{euation 2}
\end{equation}%
where $\delta _{n,k}$ is the Kronecker delta. If we take as 
\begin{equation*}
f_{L}\left( t\right) =\sum_{k=0}^{\infty }\left\langle L\mid
x^{k}\right\rangle \frac{t^{k}}{k!},
\end{equation*}
then we obtain 
\begin{equation*}
\left\langle f_{L}\left( t\right) \mid x^{n}\right\rangle =\left\langle
L\mid x^{n}\right\rangle
\end{equation*}
and so as linear functionals $L=f_{L}\left( t\right) $ (see \cite{Kim 2}, 
\cite{Dere 1}, \cite{Kim 11}, \cite{Kim 12}, \cite{Roman}). Additionally,
the map $L\rightarrow f_{L}\left( t\right) $ is a vector space isomorphism
from $\mathcal{P}^{\ast }$ onto $\mathcal{F}$. Henceforth, $\mathcal{F}$
will denote both the algebra of the formal power series in $t$ and the
vector space of all linear functionals on $\mathcal{P}$, and so an element $%
f\left( t\right) $ of $\mathcal{F}$ will be thought of as both a formal
power series and a linear functional. $\mathcal{F}$ will be called as umbral
algebra (see \cite{Kim 2}, \cite{Dere 1}, \cite{Kim 11}, \cite{Kim 12}, \cite%
{Roman}).

It is well-known that $\left\langle e^{yt}\mid x^{n}\right\rangle =y^{n}$.
Then, leads to the following%
\begin{equation*}
\left\langle e^{yt}\mid p\left( x\right) \right\rangle =p\left( y\right)
\end{equation*}
(see \cite{Kim 2}, \cite{Dere 1}, \cite{Kim 11}, \cite{Kim 12}, \cite%
{Maldonado}, \cite{Roman}). We want to note that for all $f\left( t\right) $
in $\mathcal{F}$%
\begin{equation}
f\left( t\right) =\sum_{k=0}^{\infty }\left\langle f\left( t\right) \mid
x^{k}\right\rangle \frac{t^{k}}{k!}  \label{euation 3}
\end{equation}%
and for all polynomial $p\left( x\right) $, 
\begin{equation}
p\left( x\right) =\sum_{k=0}^{\infty }\left\langle t^{k}\mid p\left(
x\right) \right\rangle \frac{x^{k}}{k!},  \label{euation 4}
\end{equation}%
(for details, see \cite{Kim 2}, \cite{Dere 1}, \cite{Kim 11}, \cite{Kim 12}, 
\cite{Roman}). The order $o\left( f\left( t\right) \right) $ of the power
series $f\left( t\right) \neq 0$ is the smallest integer $k$ for which $%
a_{k} $ does not vanish. It is considered $o\left( f\left( t\right) \right)
=\infty $ if $f\left( t\right) =0$. We see that $o\left( f\left( t\right)
g\left( t\right) \right) =o\left( f\left( t\right) \right) +o\left( g\left(
t\right) \right) $ and $o\left( f\left( t\right) +g\left( t\right) \right)
\geq \min \left\{ o\left( f\left( t\right) \right) ,o\left( g\left( t\right)
\right) \right\} $. The series $f\left( t\right) $ has a multiplicative
inverse, denoted by $f\left( t\right) ^{-1}$ or $\frac{1}{f\left( t\right) }$%
, if and only if $o\left( f\left( t\right) \right) =0$. Such series is
called an invertible series. A series $f\left( t\right) $ for which $o\left(
f\left( t\right) \right) =1$ is called a delta series (see \cite{Kim 2}, 
\cite{Dere 1}, \cite{Kim 11}, \cite{Kim 12}, \cite{Maldonado}, \cite{Roman}%
). For $f\left( t\right) ,g\left( t\right) \in \mathcal{F}$, we have $%
\left\langle f\left( t\right) g\left( t\right) \mid p\left( x\right)
\right\rangle =\left\langle f\left( t\right) \mid g\left( t\right) p\left(
x\right) \right\rangle $.

A delta series $f\left( t\right) $ has a compositional inverse $\overline{f}%
\left( t\right) $ such that $f\left( \overline{f}\left( t\right) \right) =%
\overline{f}\left( f\left( t\right) \right) =t$.

For $f\left( t\right) ,g\left( t\right) \in \mathcal{F}$ , we have $%
\left\langle f\left( t\right) g\left( t\right) \mid p\left( x\right)
\right\rangle =\left\langle f\left( t\right) \mid g\left( t\right) p\left(
x\right) \right\rangle $. By (\ref{euation 3}), we have%
\begin{equation}
p^{\left( k\right) }\left( x\right) =\frac{d^{k}p\left( x\right) }{dx^{k}}%
=\sum_{l=k}^{\infty }\frac{\left\langle t^{l}\mid p\left( x\right)
\right\rangle }{l!}l\left( l-1\right) \cdots \left( l-k+1\right) x^{l-k}%
\text{.}  \label{euation 5}
\end{equation}

Thus, notice that%
\begin{equation}
p^{\left( k\right) }\left( 0\right) =\left\langle t^{k}\mid p\left( x\right)
\right\rangle =\left\langle 1\mid p^{\left( k\right) }\left( x\right)
\right\rangle \text{.}  \label{euation 6}
\end{equation}

By (\ref{euation 5}), we have%
\begin{equation}
t^{k}p\left( x\right) =p^{\left( k\right) }\left( x\right) =\frac{%
d^{k}p\left( x\right) }{dx^{k}}\text{.}  \label{euation 7}
\end{equation}

So from the above%
\begin{equation}
e^{yt}p\left( x\right) =p\left( x+y\right) \text{.}  \label{euation 8}
\end{equation}

Let $S_{n}\left( x\right) $ be a polynomial with $\deg S_{n}\left( x\right)
=n$. Let $f\left( t\right) $ be a delta series and let $g\left( t\right) $
be an invertible series. Then there exists a unique sequence $S_{n}\left(
x\right) $ of polynomials such that $\left\langle g\left( t\right) f\left(
t\right) ^{k}\mid S_{n}\left( x\right) \right\rangle =n!\delta _{n,k}$ for
all $n,k\geq 0$. The sequence $S_{n}\left( x\right) $ is called the Sheffer
sequence for $\left( g\left( t\right) ,f\left( t\right) \right) $ or that $%
S_{n}\left( t\right) $ is Sheffer for $\left( g\left( t\right) ,f\left(
t\right) \right) $.

The Sheffer sequence for $\left( 1,f\left( t\right) \right) $ is called the
associated sequence for $f\left( t\right) $ or $S_{n}\left( x\right) $ is
associated to $f\left( t\right) $. The sheffer sequence for $\left( g\left(
t\right) ,t\right) $ is called the Appell sequence for $g\left( t\right) $
or $S_{n}\left( x\right) $ is Appell for $g\left( t\right) $.

Let $p\left( x\right) \in \mathcal{P}$. Then we have%
\begin{eqnarray}
\left\langle f\left( t\right) \mid xp\left( x\right) \right\rangle
&=&\left\langle \partial _{t}f\left( t\right) \mid p\left( x\right)
\right\rangle =\left\langle f%
{\acute{}}%
\left( t\right) \mid p\left( x\right) \right\rangle ,  \label{euation 9} \\
\left\langle e^{y}+1\mid p\left( x\right) \right\rangle &=&p\left( y\right)
+p\left( 0\right) ,\text{ (see \cite{Roman}).}  \notag
\end{eqnarray}

Let $S_{n}\left( x\right) $ be sheffer for $\left( g\left( t\right) ,f\left(
t\right) \right) $. Then%
\begin{eqnarray}
h\left( t\right) &=&\sum_{k=0}^{\infty }\frac{\left\langle h\left( t\right)
\mid S_{k}\left( x\right) \right\rangle }{k!}g\left( t\right) f\left(
t\right) ^{k},\text{ }h\left( t\right) \in \mathcal{F}  \notag \\
p\left( x\right) &=&\sum_{k=0}^{\infty }\frac{\left\langle g\left( t\right)
f\left( t\right) ^{k}\mid p\left( x\right) \right\rangle }{k!}S_{k}\left(
x\right) ,\text{ }p\left( x\right) \in \mathcal{P},  \notag \\
\frac{1}{g\left( \overline{f}\left( t\right) \right) }e^{y\overline{f}\left(
t\right) } &=&\sum_{k=0}^{\infty }S_{k}\left( y\right) \frac{t^{k}}{k!},%
\text{ for all }y\in 
\mathbb{C}
,  \label{euation 10} \\
f\left( t\right) S_{n}\left( x\right) &=&nS_{n-1}\left( x\right) \text{.} 
\notag
\end{eqnarray}

Also, it is well known in \cite{Roman} that%
\begin{equation}
\left\langle f_{1}\left( t\right) f_{2}\left( t\right) \cdots f_{m}\left(
t\right) \mid x^{n}\right\rangle =\sum \binom{n}{i_{1},\cdots ,i_{m}}%
\left\langle f_{1}\left( t\right) \mid x^{i_{1}}\right\rangle \cdots
\left\langle f_{m}\left( t\right) \mid x^{i_{m}}\right\rangle
\label{euation 33}
\end{equation}%
where $f_{1}\left( t\right) ,f_{2}\left( t\right) ,\cdots ,f_{m}\left(
t\right) \in \mathcal{F}$ and the sum is over all nonnegative integers $%
i_{1},\cdots ,i_{m}$ such that $i_{1}+\cdots +i_{m}=n$ (see \cite{Roman}).

In \cite{Dere 1} and \cite{Dere 2}, Dere and Simsek have studied
applications of umbral algebra to special functions. They gave some new
interesting links for further works of many mathematicians in Analytic
numbers theory and in modern classical umbral calculus. Kim \textit{et al}.
have given some properties of umbral calculus for Frobenius-Euler
polynomials \cite{Kim 2}, Euler polynomials \cite{Kim 11} and other special
functions \cite{Kim 12}. Also, they investigated some new applications of
umbral calculus associated with $p$-adic invariants integral on $%
\mathbb{Z}
_{p}$ in \cite{Kim 1}.

By the same motivation, we also give some applications of umbral calculus by
using extended fermionic $p$-adic $q$-integral on $%
\mathbb{Z}
_{p}$. From those applications, we derive some interesting equalities on
weighted $q$-Euler numbers, weighted $q$-Euler polynomials and weighted $q$%
-Euler polynomials of order $k$.

\section{\textbf{On the extended fermionic }$p$\textbf{-adic }$q$-\textbf{%
integrals on }$%
\mathbb{Z}
_{p}$\textbf{\ in connection with applications of umbral calculus}}

Suppose that $S_{n}\left( x\right) $ is an Appell sequence for $g\left(
t\right) $. Then, by (\ref{euation 10}), we have 
\begin{equation}
\frac{1}{g\left( t\right) }x^{n}=S_{n}\left( x\right) \text{ if and only if }%
x^{n}=g\left( t\right) S_{n}\left( x\right) \text{, }\left( n\geq 0\right) 
\text{.}  \label{euation 13}
\end{equation}

Let us contemplate as follows: 
\begin{equation*}
g_{q}\left( t\mid \zeta \right) =\frac{q\zeta e^{t}+1}{\left[ 2\right] _{q}}%
\in \mathcal{F}\text{.}
\end{equation*}

Therefore, we easily notice that $g\left( t\right) $ is an invertible
series. By (\ref{euation 13}), we have 
\begin{equation}
\sum_{n=0}^{\infty }E_{n,\zeta }^{q}\left( x\right) \frac{t^{n}}{n!}=\frac{1%
}{g_{q}\left( t\mid \zeta \right) }e^{xt}\text{.}  \label{euation 14}
\end{equation}

By (\ref{euation 14}), we have%
\begin{equation}
\frac{1}{g_{q}\left( t\mid \zeta \right) }x^{n}=E_{n,\zeta }^{q}\left(
x\right) \text{.}  \label{euation 15}
\end{equation}

Also, by (\ref{euation 10}), we have%
\begin{equation}
tE_{n,\zeta }^{q}\left( x\right) =\left( E_{n,\zeta }^{q}\left( x\right)
\right) 
{\acute{}}%
=nE_{n-1,\zeta }^{q}\left( x\right) \text{,}  \label{euation 16}
\end{equation}

By (\ref{euation 15}) and (\ref{euation 16}), we have the following
proposition.

\begin{proposition}
\label{Proposition 1}For $n\geq 0$, $E_{n,\zeta }^{q}\left( x\right) $ is an
Appell sequence for $g_{q}\left( t\mid \zeta \right) =\frac{\zeta qe^{t}+1}{%
\left[ 2\right] _{q}}$.
\end{proposition}

By (\ref{euation 11}), we derive that 
\begin{align}
\sum_{n=1}^{\infty }E_{n,\zeta }^{q}\left( x\right) \frac{t^{n}}{n!}& =\frac{%
xg_{q}\left( t\mid \zeta \right) e^{xt}-g%
{\acute{}}%
_{q}\left( t\mid \zeta \right) e^{xt}}{g\left( t\right) ^{2}}
\label{euation 17} \\
& =\sum_{n=0}^{\infty }\left( x\frac{1}{g_{q}\left( t\mid \zeta \right) }%
x^{n}-\frac{g%
{\acute{}}%
_{q}\left( t\mid \zeta \right) }{g_{q}\left( t\mid \zeta \right) }\frac{1}{%
g_{q}\left( t\mid \zeta \right) }x^{n}\right) \frac{t^{n}}{n!}  \notag
\end{align}

Because of (\ref{euation 15}) and (\ref{euation 17}), we discover the
following:%
\begin{equation*}
E_{n+1,\zeta }^{q}\left( x\right) =xE_{n,\zeta }^{q}\left( x\right) -\frac{g%
{\acute{}}%
_{q}\left( t\mid \zeta \right) }{g_{q}\left( t\mid \zeta \right) }E_{n,\zeta
}^{q}\left( x\right) \text{.}
\end{equation*}

Therefore, we get the following theorem.

\begin{theorem}
\label{Theorem 1}Let $g_{q}\left( t\mid \zeta \right) =\frac{\zeta qe^{t}+1}{%
2}\in \mathcal{F}$. Then we have for $n\geq 0:$%
\begin{equation}
E_{n+1,\zeta }^{q}\left( x\right) =\left( x-\frac{g%
{\acute{}}%
_{q}\left( t\mid \zeta \right) }{g_{q}\left( t\mid \zeta \right) }\right)
E_{n,\zeta }^{q}\left( x\right) \text{.}  \label{euation 19}
\end{equation}%
Moreover,%
\begin{equation*}
\lambda \left( -n-1,x:q:\zeta \right) =\left( x-\frac{g%
{\acute{}}%
_{q}\left( t\mid \zeta \right) }{g_{q}\left( t\mid \zeta \right) }\right)
\lambda \left( -n,x:q:\zeta \right) \text{.}
\end{equation*}%
where $g%
{\acute{}}%
_{q}\left( t\mid \zeta \right) =\frac{dg_{q}\left( t\mid \zeta \right) }{dt}$%
.
\end{theorem}

From (\ref{euation 11}), it is easy to show that%
\begin{equation*}
\sum_{n=0}^{\infty }\left( \zeta qE_{n,\zeta }^{q}\left( x+1\right)
+E_{n,\zeta }^{q}\left( x\right) \right) \frac{t^{n}}{n!}=\sum_{n=0}^{\infty
}\left( \left[ 2\right] _{q}x^{n}\right) \frac{t^{n}}{n!}\text{.}
\end{equation*}

By comparing the coefficients in the both sides of $\frac{t^{n}}{n!}$ on the
above, we develop the following:%
\begin{equation}
\zeta qE_{n,\zeta }^{q}\left( x+1\right) +E_{n,\zeta }^{q}\left( x\right) = 
\left[ 2\right] _{q}x^{n}\text{.}  \label{euation 21}
\end{equation}

From Theorem \ref{Theorem 1}, we get the following equation (\ref{euation 22}%
):%
\begin{equation}
g_{q}\left( t\mid \zeta \right) E_{n+1,\zeta }^{q}\left( x\right)
=g_{q}\left( t\mid \zeta \right) xE_{n,\zeta }^{q}\left( x\right) -g%
{\acute{}}%
_{q}\left( t\mid \zeta \right) E_{n,\zeta }^{q}\left( x\right) \text{.}
\label{euation 22}
\end{equation}

So from above%
\begin{equation*}
\left( \zeta qe^{t}+1\right) E_{n+1,\zeta }^{q}\left( x\right) =\left( \zeta
qe^{t}+1\right) xE_{n,\zeta }^{q}\left( x\right) -\zeta qe^{t}E_{n,\zeta
}^{q}\left( x\right) \text{.}
\end{equation*}

Thus, we can write the following equation:%
\begin{equation}
\zeta qE_{n+1,\zeta }^{q}\left( x+1\right) +E_{n+1,\zeta }^{q}\left(
x\right) =\zeta q\left( x+1\right) E_{n,\zeta }^{q}\left( x+1\right)
+xE_{n,\zeta }^{q}\left( x\right) -\zeta qE_{n,\zeta }^{q}\left( x+1\right) 
\text{.}  \label{equation 37}
\end{equation}

From (\ref{euation 21}) (\ref{euation 22}) and (\ref{equation 37}), we can
state following theorem.

\begin{theorem}
For $n\geq 0$, then we have%
\begin{equation}
\zeta qE_{n,\zeta }^{q}\left( x+1\right) +E_{n,\zeta }^{q}\left( x\right) = 
\left[ 2\right] _{q}x^{n}\text{.}  \label{euation 23}
\end{equation}
\end{theorem}

\begin{remark}
Assume that $S_{n}\left( x\right) $ is Sheffer sequence for $\left( g\left(
t\right) ,f\left( t\right) \right) $. Then Sheffer identity is introduced by 
\begin{equation}
S_{n}\left( x+y\right) =\sum_{k=0}^{n}\binom{n}{k}P_{k}\left( y\right)
S_{n-k}\left( x\right) =\sum_{k=0}^{n}\binom{n}{k}P_{k}\left( x\right)
S_{n-k}\left( y\right) \text{,}  \label{euation 55}
\end{equation}%
where $P_{k}\left( y\right) =S_{k}\left( y\right) g\left( t\right) $ is
associated to $f\left( t\right) $ (for details, see \cite{Kim 11}, \cite%
{Dere 1}, \cite{Dere 2}, \cite{Roman}).
\end{remark}

On account of (\ref{equation 35}) and (\ref{euation 55}), then we have%
\begin{eqnarray*}
E_{n,\zeta }^{q}\left( x+y\right) &=&\sum_{k=0}^{n}\binom{n}{k}P_{k}\left(
y\right) S_{n-k}\left( x\right) \\
&=&\sum_{k=0}^{n}\binom{n}{k}E_{n-k,\zeta }^{q}\left( y\right) x^{k}\text{.}
\end{eqnarray*}

So we have%
\begin{equation*}
E_{n,\zeta }^{q}\left( x+y\right) =\sum_{k=0}^{n}\binom{n}{k}E_{n-k,\zeta
}^{q}\left( y\right) x^{k}\text{.}
\end{equation*}

By (\ref{equation 35}), we easily see for $\alpha \left( \neq 0\right) \in 
\mathbb{C}
:$%
\begin{equation}
E_{n,\zeta }^{q}\left( \alpha x\right) =\frac{g_{q}\left( t\mid \zeta
\right) }{g_{q}\left( \frac{t}{\alpha }\mid \zeta \right) }E_{n,\zeta
}^{q}\left( x\right) \text{.}  \label{euation 56}
\end{equation}

From (\ref{euation 54}) and (\ref{euation 56}), we readily derive for $%
d\equiv 1\left( \func{mod}2\right) :$%
\begin{equation*}
\frac{g_{q}\left( t\mid \zeta \right) }{g_{q}\left( \frac{t}{d}\mid \zeta
\right) }E_{n,\zeta }^{q}\left( x\right) =\frac{d^{n}}{\left[ d\right] _{-q}}%
\sum_{j=0}^{d-1}\left( -1\right) ^{j}\zeta ^{j}q^{j}E_{n,\zeta
^{d}}^{q^{d}}\left( x+\frac{j}{d}\right) \text{.}
\end{equation*}

Let us consider the linear functional $f\left( t\right) $ that satisfies: 
\begin{equation}
\left\langle f\left( t\right) \mid p\left( x\right) \right\rangle =\int_{%
\mathbb{Z}
_{p}}\zeta ^{\xi }p\left( \xi \right) d\mu _{-q}\left( \xi \right) ,
\label{equation 38}
\end{equation}%
for all polynomials $p\left( x\right) $. From (\ref{equation 38}), we
readily see that%
\begin{equation}
f\left( t\right) =\sum_{n=0}^{\infty }\frac{\left\langle f\left( t\right)
\mid x^{n}\right\rangle }{n!}t^{n}=\sum_{n=1}^{\infty }\left( \int_{%
\mathbb{Z}
_{p}}\zeta ^{\xi }\xi ^{n}d\mu _{-q}\left( \xi \right) \right) \frac{t^{n}}{%
n!}=\int_{%
\mathbb{Z}
_{p}}\zeta ^{\xi }e^{\xi t}d\mu _{-q}\left( \xi \right) \text{.}
\label{equation 39}
\end{equation}

Thus, we have%
\begin{equation}
f\left( t\right) =\int_{%
\mathbb{Z}
_{p}}\zeta ^{\xi }e^{\xi t}d\mu _{-q}\left( \xi \right) =\frac{\left[ 2%
\right] _{q}}{\zeta qe^{t}+1}\text{.}  \label{equation 40}
\end{equation}

Therefore, by (\ref{equation 38}) and (\ref{equation 40}), we arrive at the
following theorem.

\begin{theorem}
For $n\geq 0$, then we have%
\begin{equation}
\left\langle f\left( t\right) \mid p\left( x\right) \right\rangle =\int_{%
\mathbb{Z}
_{p}}\zeta ^{\xi }p\left( \xi \right) d\mu _{-q}\left( \xi \right) \text{.}
\label{euation 24}
\end{equation}%
Also, 
\begin{equation}
\left\langle \frac{\left[ 2\right] _{q}}{\zeta qe^{t}+1}\mid p\left(
x\right) \right\rangle =\int_{%
\mathbb{Z}
_{p}}\zeta ^{\xi }p\left( \xi \right) d\mu _{-q}\left( \xi \right) .
\label{equation 41}
\end{equation}%
Obviously that%
\begin{equation}
E_{n,\zeta }^{q}=\left\langle \int_{%
\mathbb{Z}
_{p}}\zeta ^{\xi }e^{\xi t}d\mu _{-q}\left( \xi \right) \mid
x^{n}\right\rangle \text{.}  \label{equation 42}
\end{equation}
\end{theorem}

From (\ref{euation 11}) and (\ref{equation 42}), we see that%
\begin{equation}
\sum_{n=0}^{\infty }\left( \int_{%
\mathbb{Z}
_{p}}\zeta ^{\xi }\left( x+\xi \right) ^{n}d\mu _{-q}\left( \xi \right)
\right) \frac{t^{n}}{n!}=\int_{%
\mathbb{Z}
_{p}}\zeta ^{\xi }e^{\left( x+\xi \right) t}d\mu _{-q}\left( \xi \right)
=\sum_{n=0}^{\infty }\left( \int_{%
\mathbb{Z}
_{p}}\zeta ^{\xi }e^{\xi t}d\mu _{-q}\left( \xi \right) x^{n}\right) \frac{%
t^{n}}{n!}\text{.}  \label{euation 25}
\end{equation}

By (\ref{euation 11}) and (\ref{euation 26}), we see that for $n\in 
\mathbb{N}
^{\ast }$:%
\begin{equation}
E_{n,\zeta }^{q}\left( x\right) =\int_{%
\mathbb{Z}
_{p}}\zeta ^{\xi }\left( x+\xi \right) ^{n}d\mu _{-q}\left( \xi \right)
=\int_{%
\mathbb{Z}
_{p}}\zeta ^{\xi }e^{\xi t}d\mu _{-q}\left( \xi \right) x^{n}\text{.}
\label{euation 26}
\end{equation}

Consequently, we obtain the following theorem.

\begin{theorem}
For $p\left( x\right) \in \mathcal{P}$, then we have%
\begin{eqnarray}
\int_{%
\mathbb{Z}
_{p}}\zeta ^{\xi }p\left( x+\xi \right) d\mu _{-q}\left( \xi \right)
&=&\int_{%
\mathbb{Z}
_{p}}\zeta ^{\xi }e^{\xi t}d\mu _{-q}\left( \xi \right) p\left( x\right)
\label{euation 27} \\
&=&\frac{\left[ 2\right] _{q}}{\zeta qe^{t}+1}p\left( x\right) \text{.} 
\notag
\end{eqnarray}%
That is:%
\begin{equation}
E_{n,\zeta }^{q}\left( x\right) =\int_{%
\mathbb{Z}
_{p}}\zeta ^{\xi }e^{\xi t}d\mu _{-q}\left( \xi \right) x^{n}=\frac{\left[ 2%
\right] _{q}}{\zeta qe^{t}+1}x^{n}\text{.}  \label{equation 43}
\end{equation}
\end{theorem}

For $\left\vert 1-\zeta \right\vert _{p}<1$, we introduce weighted $q$-Euler
polynomials of order $k$ as follows: 
\begin{eqnarray}
\underset{\text{k-times}}{\underbrace{\int_{%
\mathbb{Z}
_{p}}\cdots \int_{%
\mathbb{Z}
_{p}}}}\zeta ^{\xi _{1}+\cdots +\xi _{k}}e^{\left( \xi _{1}+\cdots +\xi
_{k}+x\right) t}d\mu _{-q}\left( \xi _{1}\right) \cdots d\mu _{-q}\left( \xi
_{k}\right) &=&\left( \frac{\left[ 2\right] _{q}}{q\zeta e^{t}+1}\right)
^{k}e^{xt}  \label{equation 44} \\
&=&\sum_{n=0}^{\infty }E_{n,\zeta }^{\left( k\right) }\left( x\mid q\right) 
\frac{t^{n}}{n!}\text{.}  \notag
\end{eqnarray}%
where, for $x=0$, $E_{n,\zeta }^{\left( k\right) }\left( 0\mid q\right)
:=E_{n,\zeta }^{\left( k\right) }\left( q\right) $ are called weighted $q$%
-Euler numbers of order $k$.

By (\ref{equation 44}), we have 
\begin{eqnarray}
&&\underset{\text{k-times}}{\underbrace{\int_{%
\mathbb{Z}
_{p}}\cdots \int_{%
\mathbb{Z}
_{p}}}}\zeta ^{\xi _{1}+\cdots +\xi _{k}}\left( \xi _{1}+\cdots +\xi
_{k}+x\right) ^{n}d\mu _{-q}\left( \xi _{1}\right) \cdots d\mu _{-q}\left(
\xi _{k}\right)  \label{euation 28} \\
&=&\sum_{i_{1}+\cdots +i_{k}=n}\binom{n}{i_{1},\cdots ,i_{m}}\int_{%
\mathbb{Z}
_{p}}\zeta ^{\xi _{1}}\xi _{1}^{i_{1}}d\mu _{-q}\left( \xi _{1}\right)
\cdots \int_{%
\mathbb{Z}
_{p}}\zeta ^{\xi _{k}}\xi _{k}^{i_{k}}d\mu _{-q}\left( \xi _{k}\right) 
\notag \\
&=&\sum_{i_{1}+\cdots +i_{k}=n}\binom{n}{i_{1},\cdots ,i_{m}}E_{i_{1},\zeta
}^{q}\cdots E_{i_{k},\zeta }^{q}=E_{n,\zeta }^{\left( k\right) }\left( x\mid
q\right) \text{.}  \notag
\end{eqnarray}

Thanks to (\ref{equation 44}) and (\ref{euation 28}), we have 
\begin{equation}
E_{n,\zeta }^{\left( k\right) }\left( x\mid q\right) =\sum_{l=0}^{n}\binom{n%
}{l}x^{l}E_{n-l,\zeta }^{\left( k\right) }\left( q\right) \text{.}
\label{equation 45}
\end{equation}

From (\ref{euation 28}) and (\ref{equation 45}), we notice that $E_{n,\zeta
}^{\left( k\right) }\left( x\mid q\right) $ is a monic polynomial of degree $%
n$ with coefficients in $%
\mathbb{Q}
$. For $k\in 
\mathbb{N}
$, let us assume that%
\begin{eqnarray}
g_{q}^{\left( k\right) }\left( t\mid \zeta \right) &=&\left( \underset{\text{%
k-times}}{\underbrace{\int_{%
\mathbb{Z}
_{p}}\cdots \int_{%
\mathbb{Z}
_{p}}}}\zeta ^{\xi _{1}+\cdots +\xi _{k}}e^{\left( \xi _{1}+\cdots +\xi
_{k}\right) t}d\mu _{-q}\left( \xi _{1}\right) \cdots d\mu _{-q}\left( \xi
_{k}\right) \right) ^{-1}  \label{euation 29} \\
&=&\left( \frac{\zeta qe^{t}+1}{\left[ 2\right] _{q}}\right) ^{k}\text{.} 
\notag
\end{eqnarray}

From (\ref{euation 29}), we note that $g_{q}^{\left( k\right) }\left( t\mid
\zeta \right) $ is an invertible series. On account of (\ref{equation 44})
and (\ref{euation 29}), we readily derive that 
\begin{align}
\frac{1}{g_{q}^{\left( k\right) }\left( t\mid \zeta \right) }e^{xt}& =%
\underset{\text{k-times}}{\underbrace{\int_{%
\mathbb{Z}
_{p}}\cdots \int_{%
\mathbb{Z}
_{p}}}}\zeta ^{\xi _{1}+\cdots +\xi _{k}}e^{\left( \xi _{1}+\cdots +\xi
_{k}+x\right) t}d\mu _{-q}\left( \xi _{1}\right) \cdots d\mu _{-q}\left( \xi
_{k}\right)  \label{euation 30} \\
& =\sum_{n=0}^{\infty }E_{n,\zeta }^{\left( k\right) }\left( x\mid q\right) 
\frac{t^{n}}{n!}\text{.}  \notag
\end{align}

Also, we note that%
\begin{equation}
tE_{n,\zeta }^{\left( k\right) }\left( x\mid q\right) =nE_{n-1,\zeta
}^{\left( k\right) }\left( x\mid q\right) \text{.}  \label{equation 46}
\end{equation}

By (\ref{euation 30}) and (\ref{equation 46}), we easily see that $%
E_{n,\zeta }^{\left( k\right) }\left( x\mid q\right) $ is an Appell sequence
for $g_{q}^{\left( k\right) }\left( t\mid \zeta \right) $. Then, by (\ref%
{euation 30}) and (\ref{equation 46}), we get the following theorem.

\begin{theorem}
For $p\left( x\right) \in \mathcal{P}$ and $k\in 
\mathbb{N}
$, we have 
\begin{equation}
\underset{\text{k-times}}{\underbrace{\int_{%
\mathbb{Z}
_{p}}\cdots \int_{%
\mathbb{Z}
_{p}}}}\zeta ^{\xi _{1}+\cdots +\xi _{k}}p\left( \xi _{1}+\cdots +\xi
_{k}+x\right) d\mu _{-q}\left( \xi _{1}\right) \cdots d\mu _{-q}\left( \xi
_{k}\right) =\left( \frac{\left[ 2\right] _{q}}{\zeta qe^{t}+1}\right)
^{k}p\left( x\right) \text{.}  \label{equation 31}
\end{equation}%
In the special case, the weighted $q$-Euler polynomials of degree $k$ are
derived by%
\begin{equation*}
E_{n,\zeta }^{\left( k\right) }\left( x\mid q\right) =\left( \frac{2}{\zeta
qe^{t}+1}\right) ^{k}x^{n}=\int_{%
\mathbb{Z}
_{p}}\cdots \int_{%
\mathbb{Z}
_{p}}\zeta ^{\xi _{1}+\cdots +\xi _{k}}e^{\left( \xi _{1}+\cdots +\xi
_{k}\right) t}d\mu _{-q}\left( \xi _{1}\right) \cdots d\mu _{-q}\left( \xi
_{k}\right) x^{n}
\end{equation*}%
Thus, we get%
\begin{equation*}
E_{n,\zeta }^{\left( k\right) }\left( x\mid q\right) \sim \left( \left( 
\frac{\zeta qe^{t}+1}{\left[ 2\right] _{q}}\right) ^{k},t\right) \text{.}
\end{equation*}
\end{theorem}

Let us take the linear functional $f^{\left( k\right) }\left( t\right) $
that satisfies%
\begin{equation}
\left\langle f^{\left( k\right) }\left( t\right) \mid p\left( x\right)
\right\rangle =\int_{%
\mathbb{Z}
_{p}}\cdots \int_{%
\mathbb{Z}
_{p}}\zeta ^{\xi _{1}+\cdots +\xi _{k}}p\left( \xi _{1}+\cdots +\xi
_{k}\right) d\mu _{-q}\left( \xi _{1}\right) \cdots d\mu _{-q}\left( \xi
_{k}\right) \text{,}  \label{euation 32}
\end{equation}%
for all polynomials $p\left( x\right) $. Therefore, we compute as follows:%
\begin{eqnarray*}
f^{\left( k\right) }\left( t\right) &=&\sum_{n=0}^{\infty }\frac{%
\left\langle f^{\left( k\right) }\left( t\right) \mid x^{n}\right\rangle }{n!%
}t^{n} \\
&=&\sum_{n=0}^{\infty }\left( \int_{%
\mathbb{Z}
_{p}}\cdots \int_{%
\mathbb{Z}
_{p}}\zeta ^{\xi _{1}+\cdots +\xi _{k}}\left( \xi _{1}+\cdots +\xi
_{k}\right) ^{n}d\mu _{-q}\left( \xi _{1}\right) \cdots d\mu _{-q}\left( \xi
_{k}\right) \right) \frac{t^{n}}{n!} \\
&=&\int_{%
\mathbb{Z}
_{p}}\cdots \int_{%
\mathbb{Z}
_{p}}\zeta ^{\xi _{1}+\cdots +\xi _{k}}e^{\left( \xi _{1}+\cdots +\xi
_{k}\right) t}d\mu _{-q}\left( \xi _{1}\right) \cdots d\mu _{-q}\left( \xi
_{k}\right) \\
&=&\left( \frac{\left[ 2\right] _{q}}{\zeta qe^{t}+1}\right) ^{k}\text{.}
\end{eqnarray*}

Therefore, the following theorem can be expressed.

\begin{theorem}
For $p\left( x\right) \in \mathcal{P}$, we have 
\begin{eqnarray*}
&&\left\langle \int_{%
\mathbb{Z}
_{p}}\cdots \int_{%
\mathbb{Z}
_{p}}\zeta ^{\xi _{1}+\cdots +\xi _{k}}e^{\left( \xi _{1}+\cdots +\xi
_{k}\right) t}d\mu _{-q}\left( \xi _{1}\right) \cdots d\mu _{-q}\left( \xi
_{k}\right) \mid p\left( x\right) \right\rangle \\
&=&\int_{%
\mathbb{Z}
_{p}}\cdots \int_{%
\mathbb{Z}
_{p}}\zeta ^{\xi _{1}+\cdots +\xi _{k}}p\left( \xi _{1}+\cdots +\xi
_{k}\right) d\mu _{-q}\left( \xi _{1}\right) \cdots d\mu _{-q}\left( \xi
_{k}\right)
\end{eqnarray*}%
Furthermore,%
\begin{equation*}
\left\langle \left( \frac{\left[ 2\right] _{q}}{\zeta qe^{t}+1}\right)
^{k}\mid p\left( x\right) \right\rangle =\int_{%
\mathbb{Z}
_{p}}\cdots \int_{%
\mathbb{Z}
_{p}}\zeta ^{\xi _{1}+\cdots +\xi _{k}}p\left( \xi _{1}+\cdots +\xi
_{k}\right) d\mu _{-q}\left( \xi _{1}\right) \cdots d\mu _{-q}\left( \xi
_{k}\right) \text{.}
\end{equation*}%
That is:%
\begin{equation*}
E_{n,\zeta }^{\left( k\right) }\left( q\right) =\left\langle \int_{%
\mathbb{Z}
_{p}}\cdots \int_{%
\mathbb{Z}
_{p}}\zeta ^{\xi _{1}+\cdots +\xi _{k}}e^{\left( \xi _{1}+\cdots +\xi
_{k}\right) t}d\mu _{-q}\left( \xi _{1}\right) \cdots d\mu _{-q}\left( \xi
_{k}\right) \mid x^{n}\right\rangle \text{.}
\end{equation*}
\end{theorem}

From (\ref{euation 33}), we notice that%
\begin{eqnarray*}
&&\left\langle \int_{%
\mathbb{Z}
_{p}}\cdots \int_{%
\mathbb{Z}
_{p}}\zeta ^{\xi _{1}+\cdots +\xi _{k}}e^{\left( \xi _{1}+\cdots +\xi
_{k}\right) t}d\mu _{-q}\left( \xi _{1}\right) \cdots d\mu _{-q}\left( \xi
_{k}\right) \mid x^{n}\right\rangle \\
&=&\sum_{i_{1}+\cdots +i_{k}=n}\binom{n}{i_{1},\cdots ,i_{m}}\left\langle
\int_{%
\mathbb{Z}
_{p}}\zeta ^{\xi _{1}}e^{\xi _{1}t}d\mu _{-q}\left( \xi _{1}\right) \mid
x^{i_{1}}\right\rangle \cdots \left\langle \int_{%
\mathbb{Z}
_{p}}\zeta ^{\xi _{k}}e^{\xi _{k}t}d\mu _{-q}\left( \xi _{k}\right) \mid
x^{i_{k}}\right\rangle
\end{eqnarray*}

Therefore, we have%
\begin{equation*}
E_{n,\zeta }^{\left( k\right) }\left( q\right) =\sum_{i_{1}+\cdots +i_{k}=n}%
\binom{n}{i_{1},\cdots ,i_{m}}E_{i_{1},\zeta }^{q}\cdots E_{i_{k},\zeta }^{q}%
\text{.}
\end{equation*}

\begin{remark}
Our applications for weighted Euler polynomials, weighted $q$-Euler numbers
and weighted $q$-Euler polynomials of order $k$ seem to be interesting for
evaluating at $q=\zeta =1$ which lead to Euler polynomials and Euler
polynomials of order $k$, are defined respectively by%
\begin{eqnarray*}
\sum_{n=0}^{\infty }E_{n}\left( x\right) \frac{t^{n}}{n!} &=&\frac{2}{e^{t}+1%
}e^{xt}\text{,} \\
\sum_{n=0}^{\infty }E_{n}^{\left( k\right) }\left( x\right) \frac{t^{n}}{n!}
&=&\left( \frac{2}{e^{t}+1}\right) ^{k}e^{xt}\text{.}
\end{eqnarray*}%
Also, it is well known that they have representations in terms of fermionic $%
p$-adic integral on $%
\mathbb{Z}
_{p}$ as follows:%
\begin{eqnarray*}
E_{n}\left( x\right) &=&\int_{%
\mathbb{Z}
_{p}}\left( x+\xi \right) ^{n}d\mu _{-1}\left( \xi \right) , \\
E_{n}^{\left( k\right) }\left( x\right) &=&\int_{%
\mathbb{Z}
_{p}}\cdots \int_{%
\mathbb{Z}
_{p}}\left( \xi _{1}+\cdots +\xi _{k}+x\right) ^{n}d\mu _{-1}\left( \xi
_{1}\right) \cdots d\mu _{-1}\left( \xi _{k}\right) \text{.}
\end{eqnarray*}
\end{remark}

\end{document}